\documentclass[11pt]{article}
\usepackage{latexsym,amssymb}
\usepackage{amsmath}
\usepackage[mathscr]{eucal}
\usepackage{color}

\newtheorem{Theorem}{\bf Theorem}[section]
\newtheorem{Lemma}{\bf Lemma}[section]
\newtheorem{Proposition}{\bf Proposition}[section]  
\newtheorem{Corollary}{\bf Corollary}[section]
\newtheorem{Remark}{\bf Remark}[section]
\newtheorem{Example}{\bf Example}[section]
\newtheorem{Definition}{\bf Definition}[section]

\newenvironment{theorem}{\begin{Theorem}$\!\!\!$}{\end{Theorem}}
\newenvironment{lemma}{\begin{Lemma}$\!\!\!$}{\end{Lemma}}

\newenvironment{remark}{\begin{Remark}$\!\!\!$}{\end{Remark}}

\numberwithin{equation}{section}
\title{Large time behavior of ODE type solutions\\ to nonlinear diffusion equations}
\author{Junyong Eom and Kazuhiro Ishige}
\date{}
\begin{document}
\maketitle
\begin{abstract}
Consider the Cauchy problem for a nonlinear diffusion equation
\begin{equation}
\tag{P}
\left\{
\begin{array}{ll}
\partial_t u=\Delta u^m+u^\alpha & \quad\mbox{in}\quad{\bf R}^N\times(0,\infty),\vspace{5pt}\\
u(x,0)=\lambda+\varphi(x)>0 & \quad\mbox{in}\quad{\bf R}^N,\vspace{3pt}
\end{array}
\right.
\end{equation}
where $m>0$, $\alpha\in(-\infty,1)$, $\lambda>0$ and 
$\varphi\in BC({\bf R}^N)\,\cap\, L^r({\bf R}^N)$ with $1\le r<\infty$ and $\inf_{x\in{\bf R}^N}\varphi(x)>-\lambda$. 
Then the positive solution to problem~(P) behaves like a positive solution to  
ODE $\zeta'=\zeta^\alpha$ in $(0,\infty)$ and it tends to $+\infty$ as $t\to\infty$. 
In this paper we obtain the precise description of the large time behavior of the solution 
and reveal the relationship between the behavior of the solution and the diffusion effect 
the nonlinear diffusion equation has. 
\end{abstract}
\vspace{40pt}
\noindent Addresses:

\smallskip
\noindent 
J. E.: Mathematical Institute, Tohoku University,
Aoba, Sendai 980-8578, Japan. \\
\noindent 
E-mail: {\tt eom.junyong.r2@dc.tohoku.ac.jp}\\

\smallskip
\noindent K. I.:  
Graduate School of Mathematical Sciences, The University of Tokyo, 3-8-1 
Komaba, Meguro-ku, Tokyo 153-8914, Japan.\\
\noindent 
E-mail: {\tt ishige@ms.u-tokyo.ac.jp}\\

\newpage
%%%%%%%%%%%%%%%%%%%%%%%%%%%%%%
%%%%%%%%%%%%%%%%%%%%%%%%%%%%%%
\section{Introduction}
%%%%%%%%%%%%%%%%%%%%%%%%%%%%%%
%%%%%%%%%%%%%%%%%%%%%%%%%%%%%%
Let $u$ be a positive solution to the Cauchy problem 
for a nonlinear diffusion equation
\begin{equation}
\tag{P}
\left\{
\begin{array}{ll}
\partial_t u=\Delta u^m+u^\alpha\quad & \mbox{in}\quad{\bf R}^N\times(0,\infty),\vspace{5pt}\\
u(x,0)=\lambda+\varphi(x)>0 & \mbox{in}\quad{\bf R}^N,
\end{array}
\right.
\end{equation}
where $m>0$, $\alpha\in(-\infty,1)$, $\lambda>0$ and 
$\varphi\in BC({\bf R}^N)\cap L^r({\bf R}^N)$ with $1\le r<\infty$ and
$\inf_{x\in{\bf R}^N}\varphi>-\lambda$. 
For $\mu>0$, 
let $\zeta_\mu=\zeta_\mu(t)$ satisfy
$\zeta'=\zeta^\alpha$ in $(0,\infty)$ and $\zeta(0)=\mu$, 
that is, 
\begin{equation}
\label{eq:1.1}
\zeta_\mu(t)=\left(\mu^{1-\alpha}+(1-\alpha)t\right)^{\frac{1}{1-\alpha}}.
\end{equation}
Then the solution~$u$ to problem~(P) satisfies $\lim_{t\to\infty}u(x,t)=\infty$ uniformly for $x\in{\bf R}^N$ and 
\begin{equation}
\label{eq:1.2}
\lim_{t\to\infty}\,\sup_{x\in{\bf R}^N}\,\left|\frac{u(x,t)}{\zeta_\lambda(t)}-1\right|=0.
\end{equation}
Indeed, 
the comparison principle implies that 
\begin{equation}
\label{eq:1.3}
\zeta_c(t)\le u(x,t)\le\zeta_{c'}(t)\quad\mbox{in}\quad{\bf R}^N\times(0,\infty),
\end{equation}
where $c:=\inf_{x\in{\bf R}^N}u(x,0)>0$ and $c':=\sup_{x\in{\bf R}^N}u(x,0)<\infty$. 
Furthermore, by \eqref{eq:1.1} we have
$$
u(x,t)\thicksim (1-\alpha)^{\frac{1}{1-\alpha}}t^{\frac{1}{1-\alpha}}
\thicksim\zeta_\lambda(t)\quad\mbox{as}\quad t\to\infty
$$
uniformly for  $x\in{\bf R}^N$, and \eqref{eq:1.2} holds. 
(See also \cite{AE, Kaji, U, WY} for related results.) 
In this paper we obtain the precise description of the large time behavior of the solution 
and show that the function $U=U(x,t)$ defined by 
\begin{equation}
\label{eq:1.4}
U(x,t):=\frac{\lambda^\alpha\left[u(x,t)-\zeta_\lambda(t)\right]}{\zeta_\lambda(t)^{\alpha}},
\qquad x\in{\bf R}^N,\,\, t\ge 0,
\end{equation}
behaves like a rescaled solution to the heat equation as $t\to\infty$ if $m\ge\alpha$. 
\vspace{3pt}

We introduce some notation. 
For any $K\ge 0$, we denote by $[K]$ the integer such that $K-1<[K] \leq K$. 
For any $\nu  = (\nu _1,\dots,\nu _{N})\in {\bf M}:=({\bf N} \cup \{0\})^{N}$ and $x=(x_1,\dots,x_N)\in{\bf R}^N$, 
we write
$$
|\nu| := \sum_{i=1}^{N} \nu_i, \quad x^{\nu}  := x_1^{\nu_1} \cdots x_{N}^{\nu_{N}}, 
\quad \nu! := \nu_1 ! \cdots \nu_{N}!, \quad  
\partial _x^{\nu} := \frac{\partial  ^{|\nu|}}{\partial  x_1^{\nu_1} \cdots \partial  x_{N}^{\nu_{N}}}.
$$
For any $\omega=(\omega_1,\dots,\omega_N)\in{\bf M}$, 
we say $\nu\le\omega$ if $\nu_i\le\omega_i$ for $i=1,\dots,N$. 
Let $G=G(x,t)$ be the Gauss kernel, that is, 
$$
G(x,t):=(4\pi  t)^{-\frac{N}{2}} \exp\Big(-\frac{|x|^2}{4t}\Big).
$$
For any $\varphi\in BC({\bf R}^N)$, we set 
$$
\left[e^{t\Delta}\varphi\right](x):=\int_{{\bf R}^N}G(x-y,t)\varphi(y)\,dy,
\quad x\in{\bf R}^N,\,\,\,t>0. 
$$

The large time behavior of ODE type solutions was studied in \cite{IKob}, 
which deals with the Cauchy problem for the heat equation with absorption
\begin{equation}
\label{eq:1.5}
\left\{
\begin{array}{ll}
\partial_t v=\Delta v-v^\gamma\quad & \mbox{in}\quad{\bf R}^N\times(0,\infty),\vspace{5pt}\\
v(x,0)=\lambda+\varphi(x)>0 & \mbox{in}\quad{\bf R}^N,
\end{array}
\right.
\end{equation}
where $\gamma>1$, $\lambda>0$ and 
$\varphi \in BC({\bf R}^N)\cap L^r({\bf R}^N)$ with $1\le r<\infty$ and
$\inf_{x\in{\bf R}^N}\varphi(x)>-\lambda$. 
A similar argument as in \eqref{eq:1.2} implies that 
the solution~$v$ to \eqref{eq:1.5} satisfies 
$\lim_{t\to\infty}v(x,t)=0$ uniformly for $x\in{\bf R}^N$ and 
$$
\lim_{t\to\infty}\,\sup_{x\in{\bf R}^N}\,\left|\frac{v(x,t)}{\eta_\lambda(t)}-1\right|=0
\quad\mbox{with}\quad
\eta_\lambda(t):=\left(\lambda^{-(\gamma-1)}+(\gamma-1)t\right)^{-\frac{1}{\gamma-1}}
$$
(see \cite[Proposition~3.1]{IKob}). 
Here $\eta_\lambda$ is a solution of $\eta'=-\eta^\gamma$ in $(0,\infty)$ with $\eta(0)=\lambda$. 
In \cite{IKob} the second author of this paper and Kobayashi introduced a function 
$V=V(x,t)$ defined by 
$$
V(x,t):=\frac{\lambda^\gamma[v(x,t)-\eta_\lambda(t)]}{\eta_\lambda(t)^\gamma},
\qquad x\in{\bf R}^N,\,\, t\ge 0,
$$
and proved the following:
\begin{itemize}
  \item[{\rm (a)}] 
  Assume that $\varphi\in L^r({\bf R}^N)$ with $1<r<\infty$. Then 
  $$
  t^{\frac{N}{2}\left(\frac{1}{r}-\frac{1}{q}\right)}\left\|V(t)-e^{t\Delta}\varphi\right\|_{L^q({\bf R}^N)}
  =O\left(t^{-\frac{N}{2r}}\right) + O\left(t^{-\frac{N}{2}(1-\frac{1}{r})}\right) 
  $$
  as $t\to\infty$, for any $q\in[r,\infty]$;
  \item[{\rm (b)}] 
  Assume that $\varphi\in L^1({\bf R}^N,(1+|x|^K)\,dx)$ for some $0\le K<N+2$. 
  Then there exist constants $\{M_\nu\}_{|\nu|\le K}$ such that 
  $$
  t^{\frac{N}{2}\left(1-\frac{1}{q}\right)}\biggr\|\,V(t)-\sum_{|\nu|\le K}M_\nu \partial_x^\nu G(t)\,\biggr\|_{L^q({\bf R}^N)}
  =\left\{
  \begin{array}{ll}
  O\left(t^{-\frac{K}{2}}\right) & \quad\mbox{if}\quad K>[K],\vspace{5pt}\\
  o\left(t^{-\frac{K}{2}}\right) & \quad\mbox{if}\quad K=[K],
  \end{array}
  \right.
  $$
  as $t\to\infty$, for any $q\in[1,\infty]$. 
\end{itemize}
We can improve statement~(b) with the aid of \cite{Ishige14} 
and we obtain:  
\begin{itemize}
  \item[{\rm (b')}] 
  Assume that $\varphi\in L^1({\bf R}^N,(1+|x|^K)\,dx)$ for some $0\le K<N+2$. 
  Then there exist constants $\{M_\nu\}_{|\nu|\le K}$ such that 
  $$
  t^{\frac{N}{2}\left(1-\frac{1}{q}\right)}\biggr\|\,V(t)-\sum_{|\nu|\le K}M_\nu \partial_x^\nu G(t)\,\biggr\|_{L^q({\bf R}^N)}
  =o\left(t^{-\frac{K}{2}}\right)
  $$
  as $t\to\infty$, for any $q\in[1,\infty]$. 
\end{itemize}
\noindent
By statements~(a), (b) and (b') 
we see that the function~$V$ behaves like a solution of the heat equation 
as $t\to\infty$. 
In particular, by statements~(b) and (b') 
we obtain the higher order asymptotic expansion of the function $V$ as $t\to\infty$.  

In this paper we say that $u$ is a solution to problem~(P) if  
\begin{equation*}
\begin{split}
 & u\in C({\bf R}^N\times[0,\infty))\,\cap\, C^{2,1}({\bf R}^N\times(0,\infty)),\\
 & 0<\inf_{x\in{\bf R}^N}u(x,t)\le\sup_{x\in{\bf R}^N}u(x,t)<\infty\quad\mbox{for any $t>0$}.
\end{split}
\end{equation*}
The main purpose of this paper is to obtain the precise description of the large time behavior of 
the function $U$ (see \eqref{eq:1.4}). 
There are many results on the large time behavior of the solutions to the porous medium equation 
$\partial_t u=\Delta u^m$ (see e.g. \cite{V2}), however
it still seems difficult to obtain the precise description of the large time behavior of the solutions 
such as the higher order asymptotic expansions of the solutions. 
In this paper, for the case of $m\ge\alpha$, 
taking an advantage of the fact that the solution behaves like a solution to ODE 
$\zeta'=\zeta^{\alpha}$ 
and developing the arguments in \cite{Ishige13, Ishige14, IKob}, 
we show that the function~$U$ behaves like a rescaled solution  to the heat equation as $t\to\infty$
and obtain the higher order asymptotic expansions of $U$.
\vspace{5pt}

We state the main results of this paper. 
In Theorem~\ref{Theorem:1.1} we study the large time behavior of $U$ 
in the case of $\varphi\in L^r({\bf R}^N)$ with $r>1$ 
and we show that 
$U$ behaves like a rescaled function of $e^{t\Delta}\varphi$ as $t\to\infty$ 
if $m\ge\alpha$. 
\begin{theorem}
\label{Theorem:1.1}
Let $u$ be a solution to {\rm (P)}, where $m>0$ and $\alpha\in(-\infty,1)$ with $m \ge \alpha$. 
Let $U$ be as in \eqref{eq:1.4} and set 
\begin{equation} 
\label{eq:1.6}
\sigma(t):=\int_0^t  m\zeta_\lambda(s)^{m-1} \,ds,\quad t\ge 0.
\end{equation} 
Assume $\varphi\in BC({\bf R}^N) \cap L^r({\bf R}^N)$ for some $r\in(1,\infty)$. 
Then
\begin{equation} 
\label{eq:1.7}
\sigma(t)^{\frac{N}{2}(\frac{1}{r}-\frac{1}{q})}
\big\|\,U(t)-e^{\sigma(t) \Delta}\varphi\,\big\|_{L^q({\bf R}^N)}
=O\left(\sigma(t)^{-\frac{N}{2r}}\right)+O\left(\sigma(t)^{-\frac{N}{2}(1-\frac{1}{r})}\right)
\end{equation}
as $t\to\infty$, for any $q\in[r,\infty]$.
\end{theorem}
In Theorem~\ref{Theorem:1.2},  
under the assumption that $m\ge\alpha$ and $\varphi\in L^1({\bf R}^N,(1+|x|^K)\,dx)$
with $K\ge 0$,  
we prove that $U$ behaves like a suitable multiple of a rescaled Gauss kernel as $t\to\infty$ 
and obtain the higher order asymptotic expansions of $U$.
\begin{theorem} 
\label{Theorem:1.2}
Let $u$ be a solution to {\rm (P)}, where $m>0$ and $\alpha\in(-\infty,1)$ with $m \ge \alpha$. 
Assume $\varphi \in BC({\bf R}^N) \cap L^1({\bf R}^N, (1+|x|)^K)$ for some $K \geq 0$ with  
\begin{equation}
\label{eq:1.8} 
0 \le K < N + 2\frac{1- \alpha}{m - \alpha}\quad\mbox{if}\quad m>\alpha
\qquad\mbox{and}\qquad
0\le K<\infty\quad\mbox{if}\quad m=\alpha.
\end{equation}
Let $U$ and $\sigma$ be as in \eqref{eq:1.4} and \eqref{eq:1.6}, respectively. 
Then there exist constants $\{M_\nu\}_{|\nu|\le K}$ such that
\begin{equation} 
\label{eq:1.9}
\biggr\|\,U(t)-\sum_{|\nu| \leq K} M_{\nu}\partial_x^\nu G(\sigma(t))\,\biggr\|_{L^q({\bf R}^N)}
=o\left(\sigma(t)^{-\frac{N}{2}\left(1-\frac{1}{q}\right)-\frac{K}{2}}\right)
\end{equation}
as $t\to\infty$, for any $q\in[1,\infty]$.
\end{theorem}
In Theorem~\ref{Theorem:1.2} 
the asymptotic profile of $U$ is described by the Gauss kernel $G$ and its derivatives 
and it is independent of the Barenblatt solutions even if $m\not=1$.  
\begin{remark}
\label{Remark:1.1}
{\rm (i)} Assume that $m\ge \alpha$. 
Then it follows from \eqref{eq:1.1} that 
$$
\int_0^\infty \zeta_\lambda(s)^{m-1}\,ds=\infty.
$$ 
This implies that $\sigma(t)\to\infty$ as $t\to\infty$. 
\vspace{3pt}
\newline
{\rm (ii)} For the case of $m<\alpha$, 
$U$ converges to a nontrivial function as $t\to\infty$ 
and it does not necessarily decay as $t\to\infty$. 
See Remark~{\rm\ref{Remark:5.1}}.  
\end{remark}
Although Theorems~\ref{Theorem:1.1} and \ref{Theorem:1.2} are new even in the case of $m=1$, 
the novelty of this paper is to obtain the precise description of the large time behavior of $U$ in the case of $m\not=1$. 
For the solution $u$ to (P), we set $w(x,\tau):=U(x,t(\tau))$, 
where $t(\tau)$ is the inverse function of $\sigma(t)$ (see \eqref{eq:1.6}), that is,
$$
\int_0^{t(\tau)} m\zeta_\lambda(s)^{m-1}\,ds=\tau,\quad \tau\ge 0.
$$
Then $w$ satisfies 
\begin{equation} 
\label{eq:1.10}
\left\{
\begin{array}{ll}
\partial _{\tau}w=\mbox{div}\,(A(\tau,w)\nabla w)+F(\tau, w)\vspace{5pt}\\
\qquad=\Delta w+\mbox{div}\,H(x,\tau,w,\nabla w)+F(\tau, w)
 & \quad\mbox{in}\quad{{\bf R}}^N\times(0, \infty),\vspace{5pt}\\
 w(x,0)=\varphi(x) & \quad\mbox{in}\quad{\bf R}^N,
\end{array}
\right.
\end{equation}
where
\begin{equation}
\label{eq:1.11}
\begin{split}
 & A(\tau,w):=\Big(1+\lambda^{-\alpha}\eta(\tau)^{\alpha-1}w(x,\tau)\Big)^{m-1},\\
 & F(\tau,w):=\lambda^{\alpha} m\eta(\tau)^{-(m-1)}
\Big(( 1+ \lambda^{-\alpha} \eta(\tau)^{\alpha-1}w)^{\alpha}-1-\lambda^{-\alpha}\alpha\eta(\tau)^{\alpha-1}w\Big),\\
 & H(\tau,w,\nabla w):=(A(\tau,w)-1)\nabla w,\\
 & \eta(\tau):=\zeta_\lambda(t(\tau)).
\end{split}
\end{equation}
Here we remark that 
$H(\tau,w,\nabla w)=0$ if $m=1$ and
\begin{equation}
\label{eq:1.12}
1+\lambda^{-\alpha}\eta(\tau)^{\alpha-1}w(x,\tau)
=\frac{u(x,t(\tau))}{\zeta_\lambda(t(\tau))}.
\end{equation}
In \cite{Ishige13}, developing the arguments in \cite{Ishige09} and \cite{Ishige12},  
the second author of this paper and Kawakami 
studied the Cauchy problem to nonlinear diffusion equations of the form 
$$
\partial_\tau v=\Delta v+f(x,\tau,v,\nabla v)\quad\mbox{in}\quad{\bf R}^N\times(0,\infty)
$$
and established a method to obtain the higher order expansions of the solution 
behaving like a suitable multiple of the Gauss kernel as $t\to\infty$. 
(See also \cite{Ishige14} and \cite{IKM}.)
However, in the case of $m\not=1$, due to the nonlinear term $\mbox{div}\,H(\tau,w,\nabla w)$, 
we can not apply the arguments in \cite{Ishige09}--\cite{Ishige14} and \cite{IKM} to problem~\eqref{eq:1.10} directly. 
Indeed, it is difficult to apply their arguments to the Cauchy problem for nonlinear diffusion equations of the form 
$$
\partial_\tau v=\Delta v+\mbox{div}\,{\mathcal F}(x,\tau,v,\nabla v)\quad\mbox{in}\quad{\bf R}^N\times(0,\infty),
$$
since $\mbox{div}\,{\mathcal F}(x,t,v,\nabla v)$ includes the second derivatives of $v$. 
On the other hand, 
showing that $A(\tau,w)\to 1$ as $\tau\to\infty$ uniformly on ${\bf R}^N$,  
we can apply the parabolic regularity theorems to $w$
and obtain gradient estimates of~$w$. 
Then we can regard Cauchy problem~\eqref{eq:1.10} 
as the Cauchy problem to a nonlinear heat equation of the form
$$
\partial_\tau w=\Delta w+\mbox{div}\,F_1(x,\tau,w)+F_2(\tau,w)\quad\mbox{in}\quad{\bf R}^N\times(0,\infty).
$$
Furthermore, using nice properties of $H(\tau,w,\nabla w)$ (see e.g. \eqref{eq:4.1} and \eqref{eq:4.11}), 
we apply the arguments in \cite{Ishige14} and \cite{IKob} 
to obtain the precise description of the large time behavior of $w$. 
Then we can complete the proofs of Theorems~\ref{Theorem:1.1} and \ref{Theorem:1.2}. 

Our arguments are completely different from 
the arguments well-used in the study of the large time behavior of the solutions 
to nonlinear diffusion equations related to the porous medium equation. 
Compare with e.g. \cite{CA, FK, KP1, KP3, Kawa, PZ, S, V1, V2}. 
\vspace{5pt}

The rest of this paper is organized as follows. 
In Section~2 we recall some results on the higher order asymptotic expansion 
of the solutions to the heat equation with an inhomogeneous term. 
In Section~3 we obtain some preliminary lemmas on the behavior of $w$. 
In Sections~4 and 5 we prove Theorems~\ref{Theorem:1.1} and \ref{Theorem:1.2}, respectively. 
%%%%%%%%%%%%%%%%%%%%%%%%%%%%%%
%%%%%%%%%%%%%%%%%%%%%%%%%%%%%%
\section{Preliminaries}
%%%%%%%%%%%%%%%%%%%%%%%%%%%%%%
%%%%%%%%%%%%%%%%%%%%%%%%%%%%%%
In this section we introduce some notation and recall some properties of the Gauss kernel. 
Furthermore, we recall some lemmas on 
the higher order asymptotic expansion of the solutions to the heat equation with an inhomogeneous term. 
\vspace{3pt}

For any $1 \le r \le \infty$, let $\| \cdot \|_r$ be the usual norm of $L^r := L^r({\bf R}^N)$.
For any $K \ge 0$, we denote by $||| \cdot |||_K$ 
the norm of $L^1_K := L^1({\bf R}^N, (1 + |x|^K)\,dx)$, that is,   
$$
||| f |||_K := \int_{{\bf R}^N} |f(x)|(1+|x|^K)\,dx, \qquad f \in L^1_K.
$$ 
For any nonnegative functions $f$ and $g$ in $(a,\infty)$, where $a\in{\bf R}$, 
we say that $f\thicksim g$ as $t\to\infty$ if 
$\lim_{t\to\infty}f(t)/g(t)=1$. 
Furthermore, by the letter $C$
we denote generic positive constants and they may have different values also within the same line. 
\vspace{3pt}

We collect some estimates on the Gauss kernel. 
By the explicit representation of $G$ (see \eqref{eq:1.5}), 
for any $\nu \in {\bf M}$ and $j=0,1,2,\cdots$, 
we have
$$
|\partial _t^j \partial _x^{\nu} G(x, t)| \leq C t^{-\frac{N+|\nu|+2j}{2}} \Big[ 1 + \Big( \frac{|x|}{t^{1/2}} \Big)^{|\nu| + 2j} \Big]
\exp \Big( -\frac{|x|^2}{4t} \Big) 
$$
for $(x, t) \in {\bf R} ^{N} \times (0, \infty)$. 
Set
\begin{equation}
\label{eq:2.1}
g_\nu(x,t):=\frac{(-1)^{|\nu|}}{\nu!}(\partial_x^\nu G)(x,t+1). 
\end{equation}
Then it follows that
\begin{equation} 
\label{eq:2.2}
\begin{split}
 & \|g_\nu(t)\|_q\le C(1+t)^{-\frac{N}{2}\left(1-\frac{1}{q}\right)-\frac{|\nu|}{2}},\quad t>0,\\ 
 & \int_{{\bf R}^N} |x|^\ell |g_{\nu}(x,t)|\,dx\le C(1+t)^{\frac{\ell-|\nu|}{2}},\quad t>0,
\end{split}
\end{equation}  
where $1\le q\le\infty$ and $\ell\ge 0$. 
Furthermore, the Young inequality together with \eqref{eq:2.2} implies the following properties:
\begin{itemize}
\item
Let $\nu \in{\bf M}$ and $1\leq p \leq q \leq \infty$. 
Then there exists a constant $c_{|\nu|}>0$ independent of $p$ and $q$ such that 
\begin{equation} 
\label{eq:2.3}
\|\partial_x^{\nu} e^{t \Delta} \varphi\|_{q} 
\leq c_{|\nu|} t^{-\frac{N}{2}(\frac{1}{p}-\frac{1}{q})-\frac{|\nu|}{2} }\|\varphi\|_{p},
\quad t>0,
\end{equation} 
for $\varphi \in L^p({\bf R}^N)$.
In particular, $\|e^{t\Delta}\varphi\|_q\leq\|\varphi\|_q$ $(1\le q\le\infty)$ holds for $t>0$. 
\item
Let $K \ge 0$ and $\delta > 0$. 
Then there exists a constant $C>0$ such that
\begin{equation}
\label{eq:2.4}
||| e^{t \Delta} \varphi |||_K \le (1+\delta) ||| \varphi |||_K + C(1 + t^\frac{K}{2})\| \varphi \|_{1}
\end{equation}
for $\varphi \in L^1({\bf R}^N, (1 + |x|^K))$.
\end{itemize}
Let $f\in L^1_K$ $(K\ge 0)$ and $\nu\in{\bf M}$ with $|\nu|\le K$. 
Following \cite{Ishige12}, 
we define $m_\nu(f,t)$ inductively (in $\nu$) by
\begin{equation*} 
\left\{
\begin{aligned}
m_0(f,t) &:= \int_{{\bf R}^N} f(x)\, dx,\\
m_\nu(f,t) &:= \int_{{\bf R}^N} x^\nu f(x)\, dx - \sum_{\omega\leq\nu, \omega\neq\nu} 
m_\omega(f,t) \int_{{\bf R}^N} x^\nu g_\omega(x,t) \, dx \quad \mbox{if} \quad\nu\neq 0. 
\end{aligned}
\right.
\end{equation*}
Then 
\begin{equation}
\label{eq:2.5}
\int_{{\bf R}^N}
x^\omega\biggr[f-\sum_{|\nu|\le K}m_\nu(f,t)g_\nu(x,t)\biggr]\,dx=0,\qquad t\ge 0,
\end{equation}
holds for any $\omega\in{\bf M}$ with $|\omega|\le K$ 
(see \cite[Section~2]{Ishige09} and \cite[Lemma~2.1]{Ishige14}). 
In addition, due to \eqref{eq:2.5}, we have the following two lemmas 
(see \cite[Theorems~1.1 and 1.2]{Ishige14}). 
\begin{lemma}
\label{Lemma:2.1}
Let $\varphi \in L^1_K({\bf R}^N)$ for some $K\ge 0$. 
Then 
$$ 
\lim_{t\to\infty} t^{\frac{K}{2}+\frac{N}{2}(1-\frac{1}{q})} 
\biggr\|e^{t \Delta}\varphi-\sum_{|\nu| \leq K} m_{\nu}(\varphi,0) g_{\nu}(t)\biggr\|_q = 0
$$
for any $q\in[1,\infty]$.
\end{lemma}
\begin{lemma}
\label{Lemma:2.2}
Let $K\ge 0$ and $1\le q\le\infty$. 
Let $f$ be a measurable function in ${\bf R}^N\times(0,\infty)$ such that 
$$
E_{K,q}[f](t):=(1+t)^{\frac{K}{2}}\left[t^{\frac{N}{2}(1-\frac{1}{q})}\|f(t)\|_{L^q({\bf R}^N)}+\|f(t)\|_{L^1({\bf R}^N)} \right]+|||f(t)|||_K
\in L^\infty(0,T)
$$
for any $T>0$. 
\begin{itemize}
  \item[{\rm (i)}] 
  For any $\nu\in{\bf M}$ with $|\nu|\le K$, there exists a constant $C_1>0$ such that
  $$
  |m_\nu(f(t),t)|\le C_1(1+t)^{-\frac{K-|\alpha|}{2}}E_{K,q}[f](t)
  $$
  for almost all $t>0$.
  \item[{\rm (ii)}]  
  Set
  \begin{equation*}
  \begin{split}
  R_K[f](t):=
  \int_0^t e^{(t-s)\Delta}f(s)ds-\sum_{|\nu|\le K}\left[\int_0^t m_\nu(f(s),s)ds\right]g_\nu(t).
  \end{split}
  \end{equation*}
Let $j \in \{0,1\}$ and $T_*>0$. 
Then there exists a constant $C_2>0$ such that, 
for any $\epsilon>0$ and $T\ge T_*$, 
\begin{equation}
\label{eq:2.6}
\begin{split}
 & t^{\frac{N}{2}(1-\frac{1}{q})}\|\nabla^jR_K[f](t)\|_q+t^{-\frac{\ell}{2}}|||\nabla^j R_K[f](t)|||_\ell\\
 & \qquad\quad
\le\epsilon t^{-\frac{K+j}{2}}
+C_2t^{-\frac{K}{2}}\int_T^t (t-s)^{-\frac{j}{2}}E_{K,q}[f](s)ds
\end{split}
\end{equation}
for all sufficiently large $t>0$. 
In particular, if 
$$
\int_0^\infty E_{K,q}[f](s)ds<\infty,
$$
then 
\begin{equation}
\label{eq:2.7}
\lim_{t\to\infty}
t^{\frac{K}{2}}
\biggr[t^{\frac{N}{2}\left(1-\frac{1}{q}\right)}\|R_K[f](t)\|_q+t^{-\frac{\ell}{2}}|||R_K[f](t)|||_\ell\biggr]=0.
\end{equation}
\end{itemize}
\end{lemma}
%
%%%%%%%%%%%%%%%%%%%%%%%%%%%%%%
%%%%%%%%%%%%%%%%%%%%%%%%%%%%%%
\section{Preliminary estimates of solutions}
%%%%%%%%%%%%%%%%%%%%%%%%%%%%%%
%%%%%%%%%%%%%%%%%%%%%%%%%%%%%%
Assume $m \ge \alpha$. 
Let $u$ be the solution to (P). 
Due to \eqref{eq:1.3},  
the diffusion coefficient $mu^{m-1}$ 
of the nonlinear diffusion equation in (P)
is bounded and is not degenerate in ${\bf R}^N\times(0,T)$ for any $T>0$. 
By the parabolic regularity theorems we see that $u$ is smooth in ${\bf R}^N\times(0,\infty)$ 
(see also Lemma~\ref{Lemma:3.3}). 
Let $w=w(x,\tau)$ be as in Section~1, that is,  
\begin{equation}
\label{eq:3.1} 
\begin{split}
w(x,\tau):=U(x,t(\tau)) & =\frac{\lambda^\alpha\left[u(x,t(\tau))-\zeta_\lambda(t(\tau))\right]}{\zeta_\lambda(t(\tau))^{\alpha}}\\
 & \qquad\mbox{with}\quad
\int_0^{t(\tau)} m \zeta_\lambda(s)^{m-1}\,ds = \tau
\end{split}
\end{equation} 
for $x\in{\bf R}^N$ and $\tau>0$. 
Then we have: 
\begin{lemma} 
\label{Lemma:3.1}
Assume the same conditions as in Theorem~{\rm\ref{Theorem:1.1}}. 
Then
\begin{equation}
\label{eq:3.2} 
\sup_{\tau>0}\|w(\tau)\|_{\infty}<\infty. 
\end{equation} 
\end{lemma}
{\bf Proof.}
By \eqref{eq:1.3} and \eqref{eq:1.4} we have
\begin{equation}
\label{eq:3.3}
\sup_{0<\tau<T}\|w(\tau)\|_{\infty}<\infty\quad\mbox{for any $T>0$}.
\end{equation}
On the other hand, 
for any $\mu>0$, it follows from \eqref{eq:1.1} that 
$$
\zeta_\mu(t)=(1-\alpha )^{\frac{1}{1-\alpha}} t^{\frac{1}{1-\alpha}}(1+O(t^{-1}))\quad\mbox{as}\quad t\to\infty, 
$$
which implies that 
$$
\zeta_\mu(t)-\zeta_\lambda(t)=O\left(t^{\frac{\alpha}{1-\alpha}}\right)=O(\zeta_\lambda(t)^{\alpha})
\quad\mbox{as}\quad t\to\infty.
$$
Then, by \eqref{eq:1.3} we see that
$$
\|u(t)-\zeta_\lambda(t)\|_{\infty}=O(\zeta_\lambda(t)^{\alpha})
\quad\mbox{as}\quad t\to\infty.
$$
This together with \eqref{eq:1.4} implies that 
\begin{equation}
\label{eq:3.4}
\limsup_{\tau\to\infty}\|w(\tau)\|_{\infty}
=\limsup_{t\to\infty}\|U(t)\|_{\infty}<\infty. 
\end{equation}
Therefore, by \eqref{eq:3.3} and \eqref{eq:3.4} 
we obtain \eqref{eq:3.2}, and the proof is complete.
$\Box$\vspace{5pt}

By Lemma~\ref{Lemma:3.1} we apply the Taylor theorem to 
the function $F:=F(\tau,w)$ defined by \eqref{eq:1.11}. 
Then, for any $x\in{\bf R}^N$ and $\tau>0$, 
we can find $\theta \in (0,1)$ such that
\begin{equation*}
F(\tau, w(x,\tau)) =  \frac{\alpha(\alpha-1)}{2\lambda^{\alpha}} 
m\eta (\tau)^{-(m-1)}
\Big( 1+ \theta \lambda^{-\alpha} \eta (\tau)^{\alpha-1} w(x,\tau) \Big)^{\alpha -2} \Big( \eta (\tau)^{\alpha-1} w(x,\tau) \Big)^2. 
\end{equation*}
Note that
\begin{equation}
\label{eq:3.5} 
\left\{
\begin{array}{ll}
\eta(\tau) \sim c_m\tau^{\frac{1}{m-\alpha}} & \quad \mbox{if} \quad m>\alpha,\vspace{3pt}\\
\log \eta(\tau) \sim d_m\tau & \quad \mbox{if} \quad m=\alpha, 
\end{array}
\right.
\end{equation}
as $\tau \rightarrow \infty$ for some positive constants $c_m$, $d_m$
and $\lim_{\tau\to\infty}\eta(\tau)^{\alpha-1}=0$. 
On the other hand, 
by \eqref{eq:1.3} and \eqref{eq:1.12}
we can find $C>0$ such that 
\begin{equation}
\label{eq:3.6}
\begin{split}
1\le 1+\theta\lambda^{-\alpha}\eta(\tau)^{\alpha-1}w\le 1+\lambda^{-\alpha}\eta(\tau)^{\alpha-1}w
=\frac{u(x,t(\tau))}{\zeta_\lambda(t(\tau))}\le\frac{\zeta_{c'}(t(\tau))}{\zeta_\lambda(t(\tau))}\le C & \\
\mbox{if}\quad w(x,\tau)\ge 0, & \\
1\ge 1+\theta\lambda^{-\alpha}\eta(\tau)^{\alpha-1}w\ge 1+\lambda^{-\alpha}\eta(\tau)^{\alpha-1}w
=\frac{u(x,t(\tau))}{\zeta_\lambda(t(\tau))}\ge\frac{\zeta_c(t(\tau))}{\zeta_\lambda(t(\tau))}\ge C & \\
\quad\mbox{if}\quad w(x,\tau)<0, &
\end{split}
\end{equation}
for $x\in{\bf R}^N$ and $\tau>0$.
These imply that  
\begin{equation} 
\label{eq:3.7}
\frac{F(\tau, w(x,\tau))}{w(x,\tau)^2}
\le C\eta (\tau)^{-(m-1)} \eta (\tau)^{2(\alpha-1)}
\le Ch(\tau)\quad\mbox{in}\quad{\bf R}^N\times(0,\infty),
\end{equation}
where 
\begin{equation} 
\label{eq:3.8}
h(\tau):=(1+\tau)^{-1 -\frac{1 - \alpha }{m-\alpha}}\quad\mbox{if}\quad m>\alpha,
\qquad
h(\tau):=e^{-d_m(1-\alpha)\tau}\quad\mbox{if}\quad m=\alpha.
\end{equation}
In particular, 
by Lemma~\ref{Lemma:3.1}, \eqref{eq:3.7} and \eqref{eq:3.8} 
we can find $D>1$ such that  
\begin{equation} 
\label{eq:3.9}
|F(\tau, w(x,\tau))|\le Ch(\tau)|w(x,\tau)|\le C_*(1+\tau)^{-D}|w(x,\tau)|
\quad\mbox{in}\quad{\bf R}^N\times(0,\infty),
\end{equation}
where $C_*>0$. 
Then we have the following decay estimates of $w$.   
\begin{lemma}
\label{Lemma:3.2}
Assume the same conditions as in Theorem~{\rm\ref{Theorem:1.1}}. 
Let $w$ be as in \eqref{eq:3.1}. 
  \begin{itemize}
    \item[{\rm (i)}] If $\varphi\in BC({\bf R}^N)\cap L^r({\bf R}^N)$ for some 
           $r\in[1,\infty)$, then
     \begin{equation} 
     \label{eq:3.10}
       \sup_{\tau>0}\,\tau^{\frac{N}{2}(\frac{1}{r}-\frac{1}{q})}\|w(\tau)\|_{q}<\infty
       \quad\mbox{for any $q\in[r,\infty]$}.
     \end{equation}
    \item[{\rm (ii)}] If $\varphi\in BC({\bf R}^N)\cap L^1({\bf R}^N,(1+|x|)^K\,dx)$ for some $K\ge 0$, then
     \begin{equation} 
     \label{eq:3.11}
       \sup_{\tau>0}\, (1 + \tau)^{-\frac{\ell}{2}} |||w(\tau)|||_\ell<\infty
       \quad\mbox{for any $\ell\in[0,K]$}.
     \end{equation}
\end{itemize}
\end{lemma}
{\bf Proof.}
Let $w^\pm$ satisfy
\begin{equation*}
\left\{
\begin{array}{ll}
\partial_\tau w=\mbox{div}\,(A(\tau, w(x,\tau)) \nabla w) 
\pm C_*(1+\tau)^{-D} w
 & \quad\mbox{in}\quad{\bf R}^N\times(0,\infty),\vspace{3pt}\\
w(x,0)=\max\{\pm\varphi,0\} & \quad\mbox{in}\quad{\bf R}^N,
\end{array}
\right.
\end{equation*}
where $C_*$ and $D$ are as in \eqref{eq:3.9}. 
Applying the comparison principle to \eqref{eq:1.10}, we have
\begin{equation}
\label{eq:3.12}
-w^-(x,\tau)\le w(x,\tau)\le w^+(x,\tau)
\quad\mbox{in}\quad{\bf R}^N\times(0,\infty). 
\end{equation}
Set 
$$
W^{\pm}(x,\tau) := w^{\pm}(x,\tau) \exp\left(\int_0^{\tau} \mp C_*(1+s)^{-D}\,ds\right)
\quad\mbox{in}\quad{\bf R}^N\times(0,\infty).
$$
Then $W^{\pm}$ satisfies 
\begin{equation}
\label{eq:3.13}
\partial_{\tau} W^{\pm} = \mbox{div}\,( A(\tau, w(x,\tau)) \nabla W^{\pm}) \quad\mbox{in}\quad{\bf R}^N\times(0,\infty).                        
\end{equation}
Let $\Gamma=\Gamma(x,y,t)$ be the fundamental solution to parabolic equation~\eqref{eq:3.13}. 
Since it follows from \eqref{eq:1.3}, \eqref{eq:1.11} and \eqref{eq:1.12} that 
\begin{equation}
\label{eq:3.14}
\lambda_1\le A(\tau,w(x,\tau))\le\lambda_2\quad\mbox{in}\quad{\bf R}^N\times(0,\infty)
\end{equation}
for some positive constants $\lambda_1$ and $\lambda_2$, 
we observe that 
$$
C^{-1}(\tau-s)^{-\frac{N}{2}}\exp\left(-\frac{C|x-y|^2}{\tau-s}\right)
\le\Gamma(x,\tau;y,s)\le C(\tau-s)^{-\frac{N}{2}}\exp\left(-\frac{|x-y|^2}{C(\tau-s)}\right)
$$
for $x$, $y\in{\bf R}^N$ with $x\not=y$ and $\tau>s\ge 0$. 
(See e.g. \cite[Theorem~7]{Ar}.)  
Then, combining the fact that $D>1$, we see that  
\begin{equation*}
\begin{split}
|w^\pm(x,\tau)| & \le C|W^{\pm}(x,\tau)|
\le C\int_{{\bf R}^N}\Gamma(x,\tau;y,0)|\varphi(y)|\,dy\\
 & \le C\tau^{-\frac{N}{2}}\int_{{\bf R}^N}\exp\left(-\frac{|x-y|^2}{C\tau}\right)|\varphi(y)|\,dy
\le C\left[e^{\frac{C}{4}\tau\Delta}|\varphi|\right](x)
\end{split}
\end{equation*}
for $x\in{\bf R}^N$ and $\tau>0$. 
Therefore, 
by \eqref{eq:2.3}, \eqref{eq:2.4} and \eqref{eq:3.12} 
we obtain \eqref{eq:3.10} and \eqref{eq:3.11}. 
Thus Lemma~\ref{Lemma:3.2} follows. 
$\Box$\vspace{3pt}
\newline
At the end of this section we obtain uniform estimates on $\nabla w$. 
\begin{lemma}
\label{Lemma:3.3}
Assume the same conditions as in Theorem~{\rm\ref{Theorem:1.1}}. 
Let $\varphi\in BC({\bf R}^N)\cap L^r({\bf R}^N)$ for some $r\in[1,\infty)$. 
Then 
\begin{equation}
\label{eq:3.15}
\sup_{\tau>0}\,\tau^{\frac{N}{2r}+\frac{1}{2}}\|\nabla w(\tau)\|_{\infty}<\infty.
\end{equation}
\end{lemma}
{\bf Proof.}
Let $(x_0,\tau_0)\in {\bf R}^N\times(0,\infty)$ and $L=\sqrt{\tau_0}/2$.
Set
$$
\tilde{w}(x, \tau) =L^{\frac{N}{r}} w(x_0+Lx,\tau_0+L^2\tau)
$$
for $(x,\tau)\in Q:=B(0,1)\times(-1,1)$. 
By \eqref{eq:1.10} we see that $\tilde{w}$ satisfies 
$$
\partial _{\tau} \tilde{w}-\mbox{div}\,(A(\tau_0+L^2\tau, L^{-\frac{N}{r}}\tilde{w}(x,\tau))
\nabla \tilde{w})
=L^{\frac{N}{r}+2}F(\tau_0+L^2\tau,L^{-\frac{N}{r}}\tilde{w})\quad\mbox{in}\quad Q.
$$
It follows from \eqref{eq:3.9} that
$$ 
L^{\frac{N}{r}+2}\left| F\Big(\tau_0+L^2\tau,L^{-\frac{N}{r}}\tilde{w}\Big)\right| 
\le C(1+\tau_0)^{-D}|\tilde{w} | \quad\mbox{in}\quad Q.
$$
Then Lemma~\ref{Lemma:3.2} together with \eqref{eq:3.14} implies that
$$
\sup_{\tau_0>0}\|\tilde{w}\|_{L^\infty(B(0,1)\times(-1,1))}<\infty.
$$
Then we apply the arguments in \cite[Chapter~V, Section~3]{LSU} to obtain 
$$
\sup_{\tau_0>0}\|\nabla\tilde{w}\|_{L^\infty(B(0,1/2)\times(-1/2,1/2))}<\infty,
$$
that is,  
$$
|\nabla w(x_0,\tau_0)|\le C\tau_0^{-\frac{N}{2r}-\frac{1}{2}},
\qquad x_0\in{\bf R}^N,\,\,\tau_0>0. 
$$
Therefore we have \eqref{eq:3.15}, and the proof is complete. 
$\Box$
%%%%%%%%%%%%%%%%%%%%%%%%%%%%%%
%%%%%%%%%%%%%%%%%%%%%%%%%%%%%%
\section{Proof of Theorem~\ref{Theorem:1.1}}
%%%%%%%%%%%%%%%%%%%%%%%%%%%%%%
%%%%%%%%%%%%%%%%%%%%%%%%%%%%%%
We prove Theorem~\ref{Theorem:1.1}. 
By Lemma~\ref{Lemma:3.1} and the Taylor theorem, 
for any $x\in{\bf R}^N$ and $\tau>0$, 
we can find $\theta \in (0,1)$ such that
\begin{equation*}
\begin{split}
A(\tau,w(x,\tau))-1 & =\left[1+\lambda^{-\alpha}\eta(\tau)^{\alpha-1}w(x,\tau)\right]^{m-1}-1\\
 & =(m-1)\lambda^{-\alpha}\eta(\tau)^{\alpha-1}w(x,\tau) 
\left(1+\theta \lambda^{-\alpha}\eta(\tau)^{\alpha-1}w(x,\tau)\right)^{m-1}.
\end{split} 
\end{equation*}
Then, by Lemma~\ref{Lemma:3.3} and \eqref{eq:3.6}
we have
\begin{equation}
\label{eq:4.1}
\begin{split}
|H(\tau,w(\tau),\nabla w(\tau)) |  & \le C\eta(\tau)^{\alpha-1}|w(x,\tau)||\nabla w(x,\tau)|\\
 & \le C\tau^{-\frac{N}{2r}-\frac{1}{2}}\eta(\tau)^{\alpha-1}|w(x,\tau)|
\quad\mbox{in}\quad{\bf R}^N\times(0,\infty).
\end{split}
\end{equation}
On the other hand, it follows from \eqref{eq:1.10} that
\begin{equation}
\label{eq:4.2}
w(\tau) = e^{\tau \Delta} \varphi + \int_0^{\tau} e^{(\tau-s) \Delta} F(s,w(s)) \, ds + 
\mbox{div} \, \int_0^{\tau} e^{(\tau-s) \Delta} H(s,w(s),\nabla w(s)) \,ds
\end{equation}
for $\tau>0$. 
\vspace{5pt}
\newline
{\bf Proof of Theorem~\ref{Theorem:1.1}.} 
Let $1 < r \le q \le \infty$. 
It follows from \eqref{eq:4.2} that 
\begin{equation}
\label{eq:4.3}
\begin{split}
 & \left\|w(\tau) - e^{\tau \Delta} \varphi\right\|_q\\
 & \le\left\| \int_0^{\tau} e^{(\tau-s) \Delta} F(s,w(s)) \, ds \right\|_q 
+\left\| \mbox{div} \, \int_0^{\tau} e^{(\tau-s) \Delta} H(s,w(s),\nabla w(s)) \, ds \right\|_q
\end{split}
\end{equation}
for $ \tau > 0$. 
By \eqref{eq:3.7} and \eqref{eq:3.10} we have
\begin{equation}
\label{eq:4.4}
|F(\tau,w(x,\tau))| \le Ch(\tau)|w(x,\tau)|^2 
\le C \tau^{-\frac{N}{2r}} h(\tau)|w(x,\tau)|
\quad\mbox{in}\quad{\bf R}^N\times(0,\infty).
\end{equation}
Then, by \eqref{eq:2.3}, \eqref{eq:3.8}, \eqref{eq:3.10} and \eqref{eq:4.4}
we see that
\begin{equation}
\label{eq:4.5}
\begin{split}
 & \left\| \int_{\tau/2}^{\tau} e^{(\tau-s) \Delta} F(s,w(s)) \, ds \right\|_q 
 \le \int_{\tau/2}^{\tau} \| F(s,w(s)) \|_q \, ds \\
 & \le C \int_{\tau/2}^{\tau}\tau^{ -\frac{N}{2r} }h(\tau) \| w(s) \|_q \, ds
 =O\biggr(\tau^{-\frac{N}{2}\left(\frac{2}{r} -\frac{1}{q}\right)}\biggr)
 \qquad\qquad\qquad\qquad
\end{split}
\end{equation}
for all sufficiently large $\tau>0$. 
Similarly, in the case of $r\ge 2$ we obtain 
\begin{equation}
\label{eq:4.6}
\begin{split}
 & \left\| \int_0^{\tau/2} e^{(\tau-s) \Delta} F(s,w(s)) \, ds \right\|_q
\le C \int_0^{\tau/2} (\tau-s)^{-\frac{N}{2}\left(\frac{2}{r}-\frac{1}{q}\right)} \| F(s,w(s)) \|_{r/2} \,ds \\  
 & \le C\tau^{-\frac{N}{2}\left(\frac{2}{r}-\frac{1}{q}\right)} \int_0^{\tau/2} h(s)\| w(s) \|^2_r \, ds
 =O\biggr(\tau^{-\frac{N}{2}\left(\frac{2}{r} -\frac{1}{q}\right)}\biggr)    
\end{split}
\end{equation}
for all sufficiently large $\tau>0$. 
In the case of $1<r<2$
it follows from Lemma~\ref{Lemma:3.1} and \eqref{eq:3.10} with $q=r$ that
\begin{equation}
\label{eq:4.7}
\sup_{\tau>0}\,\|w(\tau)\|_2
\le\sup_{\tau>0}\,\|w(\tau)\|_r^{\frac{r}{2}}\|w(\tau)\|_\infty^{1-\frac{r}{2}}<\infty.
\end{equation}
Then, similarly to \eqref{eq:4.6},   
\begin{equation}
\label{eq:4.8}
\begin{split}
 & \left\| \int_0^{\tau/2} e^{(\tau-s) \Delta} F(s,w(s)) \, ds \right\|_q
\le C \int_0^{\tau/2} (\tau-s)^{-\frac{N}{2}\left(1-\frac{1}{q}\right)} \| F(s,w(s)) \|_1 \,ds \\  
& \le C\tau^{-\frac{N}{2}\left(1-\frac{1}{q}\right)} \int_0^{\tau/2} h(s)\| w(s) \|_2^2 \, ds=O\biggr(\tau^{-\frac{N}{2}\left(1-\frac{1}{q}\right)}\biggr)
\end{split}
\end{equation}
for all sufficiently large $\tau>0$. 
Therefore we deduce from \eqref{eq:4.5}, \eqref{eq:4.6} and \eqref{eq:4.8} that 
\begin{equation}
\label{eq:4.9}
\tau^{\frac{N}{2}\left(\frac{1}{r}-\frac{1}{q}\right)}\left\| \int_0^{\tau} e^{(\tau-s) \Delta} F(s,w(s)) \, ds \right\|_q
=O\left(\tau^{-\frac{N}{2r}}\right)+O\left(\tau^{-\frac{N}{2}\left(1-\frac{1}{r}\right)}\right)
\quad\mbox{as}\quad\tau\to\infty.
\end{equation}
Similarly, 
by \eqref{eq:2.3}, \eqref{eq:3.5}, \eqref{eq:3.10} and \eqref{eq:4.1} we obtain
\begin{equation}
\label{eq:4.10}
\begin{split}
 & \left\| \mbox{div} \, \int_{\tau/2}^{\tau} e^{(\tau-s) \Delta} H(s,w(s),\nabla w(s)) \, ds \right\|_q 
\le  C \int_{\tau/2}^{\tau} (\tau-s)^{-\frac{1}{2}} \| H(s,w(s),\nabla w(s)) \|_q \, ds \\
 & \le C\int_{\tau/2}^{\tau} 
\tau^{-\frac{1}{2} -\frac{N}{2r}}\eta(\tau)^{\alpha-1}(\tau-s)^{-\frac{1}{2}}\| w(s) \|_q \, ds
=O\biggr(\tau^{-\frac{N}{2}\left(\frac{2}{r} -\frac{1}{q}\right)}\biggr)
\end{split}
\end{equation}
for all sufficiently large $\tau>0$. 

On the other hand, 
since 
\begin{equation}
\label{eq:4.11}
\begin{split}
H(\tau,w(\tau),\nabla w(\tau)) & =\nabla\int_0^{w(x,\tau)}
[A(\tau, \xi)-1]\,d\xi\\
 & =\nabla\int_0^{w}
[(1+\lambda^{-\alpha}\eta(\tau)^{\alpha-1} \xi )^{m-1}-1]\, d\xi,
\end{split}
\end{equation}
it follows that
\begin{equation*}
\begin{split}
&\mbox{div} \int_0^{\tau/2} e^{(\tau-s) \Delta} H(s,w(s),\nabla w(s)) \, ds \\
&=\int_0^{\tau/2}  \Delta e^{(\tau-s)\Delta}  
\left[ \int_0^{w}
[(1+\lambda^{-\alpha}\eta(\tau)^{\alpha-1} \xi )^{m-1}-1]\, d\xi \right] \, ds.
\end{split} 
\end{equation*}
Furthermore, by \eqref{eq:3.6} 
we apply the mean value theorem to obtain 
$$
\biggr|\,\int_0^{w(x,s)}
[(1+\lambda^{-\alpha}\eta(\tau)^{\alpha-1} \xi )^{m-1}-1]\, d\xi\,\biggr|
\le C\lambda^{-\alpha}\eta(\tau)^{\alpha-1}|w(x,s)|^2
$$
for $x\in{\bf R}^N$ and $s>0$. 
Then, in the case of $r\ge 2$, 
by \eqref{eq:2.3}, \eqref{eq:3.10} and \eqref{eq:4.4} we obtain
\begin{equation} 
\label{eq:4.12}
\begin{split}
&  \left\|  \mbox{div} \int_0^{\tau/2} e^{(\tau-s) \Delta} H(s,w(s),\nabla w(s)) \, ds \right\|_q  \\
&  \le C \int_0^{\tau/2} (\tau-s)^{-\frac{N}{2}\left(\frac{2}{r}-\frac{1}{q}\right)-1}
\left\| \int_0^{w}\left\{(1+\lambda^{-\alpha}\eta(\tau)^{\alpha-1} \xi )^{m-1}-1\right\}\, d\xi\right\|_{r/2} \,ds \\
& \le C\tau^{-\frac{N}{2}(\frac{2}{r}-\frac{1}{q})-1} 
\int_0^{\tau/2} \eta(\tau)^{\alpha-1}\| w(s) \|^2_r \, ds
=O\biggr(\tau^{-\frac{N}{2}\left(\frac{2}{r} -\frac{1}{q}\right)}\biggr)
\end{split}
\end{equation}
for all sufficiently large $\tau>0$. 
Similarly, in the case of $1<r<2$, 
$$
\sup_{\tau>0}\,\|w(\tau)\|_2<\infty
$$ 
(see \eqref{eq:4.7}) and we observe that 
\begin{equation}
\label{eq:4.13}
\begin{split}
& \left\|  \mbox{div} \int_0^{\tau/2} e^{(\tau-s) \Delta} H(s,w(s),\nabla w(s)) \, ds\right\|_q \\
& \le C \int_0^{\tau/2} (\tau-s)^{-\frac{N}{2}\left(1-\frac{1}{q}\right)-1}
\left\| \int_0^{w}\left\{(1+\lambda^{-\alpha}\eta(\tau)^{\alpha-1} \xi )^{m-1}-1\right\}\, d\xi\right\|_1 \,ds \\
& \le C\tau^{-\frac{N}{2}(1-\frac{1}{q}) -1}\int_0^{\tau/2} \eta(\tau)^{\alpha-1}\|w(s)\|_2^2\,ds
=O\left(\tau^{-\frac{N}{2}\left(1-\frac{1}{q}\right)}\right)
\end{split}
\end{equation}
for all sufficiently large $\tau>0$. 
Therefore we deduce from \eqref{eq:4.10}, \eqref{eq:4.12} and \eqref{eq:4.13} that 
\begin{equation}
\label{eq:4.14}
\tau^{\frac{N}{2}\left(\frac{1}{r}-\frac{1}{q}\right)}\left\| \mbox{div} \, \int_0^{\tau} e^{(\tau-s) \Delta} H(s,w(s),\nabla w(s)) \, ds \right\|_q
=O\left(\tau^{-\frac{N}{2r}}\right)+O\left(\tau^{-\frac{N}{2}\left(1-\frac{1}{r}\right)}\right)
\end{equation}
as $\tau\to\infty$. 
Consequently, by \eqref{eq:4.3}, \eqref{eq:4.9} and \eqref{eq:4.14}
we obtain 
$$
\tau^{\frac{N}{2}\left(\frac{1}{r}-\frac{1}{q}\right)}
\left\|w(\tau) - e^{\tau \Delta} \varphi\right\|_q
=O\left(\tau^{-\frac{N}{2r}}\right)+O\left(\tau^{-\frac{N}{2}\left(1-\frac{1}{r}\right)}\right)
$$
as $\tau\to\infty$. 
This implies \eqref{eq:1.7}, and the proof of Theorem~\ref{Theorem:1.1} is complete.
$\Box$
%%%%%%%%%%%%%%%%%%%%%%%%%%%%%%
%%%%%%%%%%%%%%%%%%%%%%%%%%%%%%
\section{Proof of Theorem~\ref{Theorem:1.2}}
%%%%%%%%%%%%%%%%%%%%%%%%%%%%%%
%%%%%%%%%%%%%%%%%%%%%%%%%%%%%%
We apply Lemmas~\ref{Lemma:2.1} and \ref{Lemma:2.2} to prove Theorem~\ref{Theorem:1.2}.
\vspace{3pt}
\newline
{\bf Proof of Theorem~\ref{Theorem:1.2}.}
Assume \eqref{eq:1.8}. Let $1\le q\le\infty$. 
Let $\varphi \in BC({\bf R}^N) \cap L^1_K ({\bf R}^N)$. 
Set 
\begin{equation*}
\begin{split}
 & \tilde{w}(x,\tau):=w(x,\tau+1),\\
 & \tilde{F}(x,\tau):=F(\tau+1,\tilde{w}(x,\tau)),
\quad
\tilde{H}(x,\tau):=H(\tau+1,\tilde{w}(x,\tau),\nabla\tilde{w}(x,\tau)).
\end{split}
\end{equation*}
Then it follows from \eqref{eq:4.2} that 
\begin{equation}
\label{eq:5.1}
\tilde{w}(\tau)=e^{\tau\Delta}\tilde{w}(0)
+\int_0^{\tau} e^{(\tau-s) \Delta}\tilde{F}(s)\, ds
+\,\mbox{div} \, \int_0^{\tau} e^{(\tau-s) \Delta}\tilde{H}(s)\,ds
\end{equation}
for $\tau>0$.
By Lemma~\ref{Lemma:3.2}~(ii)
we apply Lemma~\ref{Lemma:2.1} to obtain 
\begin{equation} 
\label{eq:5.2}
\lim_{\tau \to \infty} \tau^{\frac{K}{2} + \frac{N}{2}(1-\frac{1}{q})} 
\biggr\| e^{\tau \Delta}\tilde{w}(0)-\sum_{|\nu| \leq K} m_{ \nu}(\tilde{w}(0),0)g_{\nu}(\tau)\biggr\|_q = 0.
\end{equation}  
On the other hand, 
by \eqref{eq:3.7} and \eqref{eq:3.10} with $q=\infty$ we have 
$$
|\tilde{F}(x,\tau)|
\le Ch(\tau+1)|\tilde{w}(x,\tau)|^2
\le C(1+\tau)^{-\frac{N}{2}}h(\tau+1)|\tilde{w}(x,\tau)|
\quad\mbox{in}\quad{\bf R}^N\times(0,\infty).
$$ 
Then, by \eqref{eq:3.10} we see that
$$
E_{K,q}[\tilde{F}](\tau)\le C(1+\tau)^{\frac{K}{2}-\frac{N}{2}}h(\tau+1)
\quad\mbox{for $\tau>0$}. 
$$
This together with \eqref{eq:1.8} and \eqref{eq:3.8} implies that 
\begin{equation}
\label{eq:5.3}
E_{K,q}[\tilde{F}](\tau)\le C(1+\tau)^{-\gamma}\quad\mbox{for $\tau>0$ if $m>\alpha$},
\end{equation}
where 
\begin{equation}
\label{eq:5.4}
\gamma:=-\frac{K}{2}+\frac{N}{2}+\frac{1-\alpha}{m-\alpha}+1>1.
\end{equation}
Similarly, if $m=\alpha$, then \eqref{eq:5.3} holds for any $\gamma>1$. 
Then, by Lemma~\ref{Lemma:2.2}~(i), 
for any $\nu\in{\bf M}$ with $|\nu|\le K$, we have 
$$
\int_{\tau_1}^{\tau_2} |m_\nu(\tilde{F}(s),s)|\,ds
\le C\int_{\tau_1}^{\tau_2}(1+s)^{-\frac{K-|\nu|}{2}-\gamma}\,ds
\le C\tau_1^{-\frac{K-|\nu|}{2}-\gamma+1},
\quad 0<\tau_1<\tau_2.
$$
Then we can find a constant $m^{\tilde{F}}_{\nu}$ such that 
\begin{equation} 
\label{eq:5.5}
\int_0^\tau m_{\nu}(\tilde{F}(s), s)\,ds= m^{\tilde{F}}_{\nu}+O\left(\tau^{-\frac{K-|\nu|}{2}-\gamma+1}\right)
\quad\mbox{as}\quad\tau\to\infty.
\end{equation}
Furthermore, we apply Lemma~\ref{Lemma:2.2}~(ii) (see \eqref{eq:2.7}) to obtain
\begin{equation} 
\label{eq:5.6}
\biggr\| \int_0^{\tau} e^{(\tau-s) \Delta}\tilde{F}(s) \, ds -  \sum_{|{\bf \nu}| \leq K} 
\left[ \int_0^{\tau} m_{ \nu}(\tilde{F}(s) , s)\, ds \right] g_{\nu}(\tau)\biggr\|_q 
= o\left(\tau^{-\frac{K}{2}-\frac{N}{2}\left(1-\frac{1}{q}\right)}\right)
\end{equation}
as $\tau \to \infty$. 
Therefore we deduce from \eqref{eq:2.2}, \eqref{eq:5.5} and \eqref{eq:5.6} that 
\begin{equation} \label{eq:5.7}
\begin{split}
& \biggr\| \int_0^{\tau} e^{(\tau-s) \Delta} \tilde{F}(s) \, ds -  
\sum_{|{\bf \nu}| \leq K} m^{\tilde{F}}_{\nu} g_{\nu}(\tau)\biggr\|_q  \\
& \le\biggr\| \int_0^{\tau} e^{(\tau-s) \Delta}\tilde{F}(s) \, ds
 -\sum_{|{\bf \nu}| \leq K} \left[ \int_0^{\tau} m_\nu(\tilde{F}(s) , s) ds \right] g_{\nu}(\tau)\biggr\|_q \\
& \quad + \biggr\| \sum_{|{\bf \nu}| \leq K} \left[ \int_0^\tau m_\nu(\tilde{F}(s),s)\, ds 
-m^{\tilde{F}}_{\nu} \right] g_{\nu}(\tau)  \biggr\|_q 
=o\left(\tau^{-\frac{K}{2}-\frac{N}{2}\left(1-\frac{1}{q}\right)}\right)
\end{split}
\end{equation}
as $\tau \to \infty$.

Let $K':=K-1$ if $K\ge 1$ and $K':=0$ if $0<K<1$. 
By Lemma~\ref{Lemma:3.3} and \eqref{eq:4.1} we have 
$$
|\tilde{H}(x,\tau)|\le C(1+\tau)^{-\frac{N}{2}-\frac{1}{2}}\eta(\tau+1)^{\alpha-1}
|\tilde{w}(x,\tau)|
\quad\mbox{in}\quad{\bf R}^N\times(0,\infty).
$$
Then, by \eqref{eq:3.10} we obtain 
$$
E_{K',q}[\tilde{H}](\tau)\le C(1+\tau)^{\frac{K'}{2}-\frac{N}{2}-\frac{1}{2}}\eta(\tau+1)^{\alpha-1}\quad\mbox{for $\tau>0$}. 
$$
This together with \eqref{eq:1.8}, \eqref{eq:3.8} and \eqref{eq:5.4} implies that 
\begin{equation}
\label{eq:5.8}
E_{K',q}[\tilde{H}](\tau)\le C(1+\tau)^{-\gamma'}\quad\mbox{for $\tau>0$ if $m>\alpha$},
\end{equation}
where 
\begin{equation}
\begin{split}
\label{eq:5.9}
\gamma'&=-\frac{K'}{2}+\frac{N}{2}+\frac{1}{2}+\frac{1-\alpha}{m-\alpha} 
= \left\{
\begin{aligned}
&\gamma  \quad &\mbox{if}& \quad K \ge 1, \\ 
&\frac{N}{2}+\frac{1}{2}+\frac{1-\alpha}{m-\alpha}  \quad &\mbox{if}& \quad 0<K<1,
\end{aligned}
\right. \\
&> 1.
\end{split}
\end{equation}
Similarly, if $m=\alpha$, 
then \eqref{eq:5.8} holds for any $\gamma'>1$. 
Then, similarly to \eqref{eq:5.5}, 
for any $\nu\in{\bf M}$ with $|\nu|\le K'$, 
we can find a constant $m_\nu^{\tilde{H}}$ such that 
\begin{equation} 
\label{eq:5.10}
\int_0^\tau m_{\nu}(\tilde{H}(s), s)\,ds= m^{\tilde{H}}_{\nu}+O\left(\tau^{-\frac{K'-|\nu|}{2}-\gamma'+1}\right)
\quad\mbox{as}\quad\tau\to\infty.
\end{equation}
Furthermore, we apply Lemma~\ref{Lemma:2.2}~(ii) with $j=1$ (see \eqref{eq:2.6}). 
Then, 
by \eqref{eq:5.8} and \eqref{eq:5.9}, 
for any $\epsilon>0$ and any fixed sufficiently large $T>0$, 
we have 
\begin{equation}
\label{eq:5.11}
\begin{split}
 & I(\tau):=\tau^{\frac{N}{2}\left(1-\frac{1}{q}\right)}\\
 & \qquad\times\left\|\mbox{div} \, \int_0^{\tau} e^{(\tau-s) \Delta} \tilde{H}(s) \ ds
-\mbox{div} \, \sum_{|{\bf \nu}| \leq K'} \left[ \int_0^{\tau} m_\nu(\tilde{H}(s) , s)\, ds \right] g_{\nu}(\tau) \right \|_q\\
 & \qquad
\le\epsilon\tau^{-\frac{K'+1}{2}}+C\tau^{-\frac{K'}{2}}\int_T^\tau(\tau-s)^{-\frac{1}{2}}E_{K',q}[\tilde{H}](s)\,ds
\end{split}
\end{equation}
for all sufficiently large $\tau>0$. 
Since
\begin{equation*}
\begin{split}
 & \int_T^{\tau/2}(\tau-s)^{-\frac{1}{2}}E_{K',q}[\tilde{H}](s)\,ds
\le C\tau^{-\frac{1}{2}}\int_T^{\tau/2}E_{K',q}[\tilde{H}](s)\,ds,\\
 & \int_{\tau/2}^\tau(\tau-s)^{-\frac{1}{2}}E_{K',q}[\tilde{H}](s)\,ds
\le C\tau^{-\gamma'}\int_{\tau/2}^\tau(\tau-s)^{-\frac{1}{2}}\,ds\le Ct^{-\gamma'+\frac{1}{2}}
=o\left(\tau^{-\frac{1}{2}}\right), 
\end{split}
\end{equation*}
for all sufficiently large $\tau>0$, by \eqref{eq:5.11} we obtain 
\begin{equation*}
\limsup_{\tau\to\infty}\,\tau^{\frac{K}{2}}I(\tau)
\le\epsilon+C\int_T^{\tau/2} E_{K',q}[\tilde{H}](s)\,ds.
\end{equation*}
Since $\epsilon$ and $T$ are arbitrary, 
by \eqref{eq:5.8} and \eqref{eq:5.9} we observe that 
\begin{equation}
\label{eq:5.12}
\lim_{\tau\to\infty}\tau^{\frac{K}{2}}I(\tau)=0. 
\end{equation}
Therefore, by \eqref{eq:2.2}, \eqref{eq:5.10} and \eqref{eq:5.12}
we see that 
\begin{equation}
\label{eq:5.13}
\begin{split}
 & \left\|\mbox{div} \, \int_0^{\tau} e^{(\tau-s) \Delta} \tilde{H}(s) \ ds
-\mbox{div} \, \sum_{|{\bf \nu}| \leq K'}m_\nu^{\tilde{H}}g_{\nu}(\tau) \right \|_q\\
 & \le\left\|\mbox{div} \, \int_0^{\tau} e^{(\tau-s) \Delta} \tilde{H}(s) \ ds
-\mbox{div} \, \sum_{|{\bf \nu}| \leq K'} \left[ \int_0^{\tau} m_\nu(\tilde{H}(s) , s)\, ds \right] g_{\nu}(\tau) \right \|_q\\
 & \qquad
 +\left\|\mbox{div} \, \sum_{|{\bf \nu}| \leq K'} \left[ \int_0^{\tau} m_\nu(\tilde{H}(s) , s)\, ds \right] g_{\nu}(\tau)
 -\mbox{div} \, \sum_{|{\bf \nu}| \leq K'}m_\nu^{\tilde{H}}g_\nu(\tau)\right \|_q\\
 & =o\left(\tau^{-\frac{K}{2}-\frac{N}{2}\left(1-\frac{1}{q}\right)}\right)
\end{split}
\end{equation}
as $\tau\to\infty$. 
Combining \eqref{eq:2.1}, \eqref{eq:5.1}, \eqref{eq:5.2}, \eqref{eq:5.7} and \eqref{eq:5.13}, 
we can find constants $\{m_\nu\}_{|\nu|\le K}$ such that 
\begin{equation*}
\begin{split}
 & \left\|w(\tau+1)-\sum_{|\nu|\le K}\frac{(-1)^{|\nu|}}{\nu!}m_\nu 
  \partial_x^{\nu} G(\tau+1)\right\|_q\\
 & =\biggr\|\tilde{w}(\tau)-\sum_{|\nu|\le K}m_\nu g_\nu(\tau)\biggr\|_q
 =o\left(\tau^{-\frac{K}{2}-\frac{N}{2}(1-\frac{1}{q})}\right)
\end{split}
\end{equation*}
as $\tau\to\infty$. 
This implies \eqref{eq:1.9} with $M_\nu=(-1)^{|\nu|}m_\nu/\nu!$. 
Thus Theorem~\ref{Theorem:1.2} follows.
$\Box$
\begin{remark}
\label{Remark:5.1}
Assume $m<\alpha$. Then we can define $\tau^*\in(0,\infty)$ by 
$$
\tau^*:=\int_0^\infty m\zeta_\lambda(s)^{m-1}\,ds. 
$$
Then it follows from \eqref{eq:1.4} that 
$$
\lim_{t\to\infty}
\frac{\lambda^\alpha\left[u(x,t)-\zeta_\lambda(t)\right]}{\zeta_\lambda(t)^{\alpha}}
=w(x,\tau^*)\quad\mbox{for}\quad x\in{\bf R}^N.
$$
\end{remark}
\medskip

\noindent
{\bf Acknowledgements.} 
The second author of this paper was supported in part by the Grant-in-Aid for Scientific Research (A)(No.~15H02058)
from Japan Society for the Promotion of Science. 
%%%%%%%%%%%%%%%%%%%%%%%%%%%%%%
%%%%%%%%%%%%%%%%%%%%%%%%%%%%%%

%%%%%%%%%%%%%%%%%%%%%%%%%%%%%%%%%%%%%%%%%%%%%%%%%%%%%%%
\end{document}